\date{May 17, 2000} 
\title {      Values of Brownian intersection exponents III:\\
		Two-sided exponents
}
\author { Gregory F.  Lawler\thanks{Duke University, Research
supported by the National Science Foundation}
\and Oded Schramm\thanks{Microsoft Research}
\and Wendelin Werner\thanks{Universit\'e Paris-Sud}
}
\documentclass[12pt]{article}
\usepackage{amsmath}
\usepackage{amsthm}
\usepackage{amsfonts}
\newif\ifhyper\IfFileExists{hyperref.sty}{\hypertrue}{\hyperfalse}
\hyperfalse  % should be \hyperfalse when submitting to arxiv.com
\ifhyper\usepackage[naturalnames]{hyperref}\fi

\newif\ifdraft
\drafttrue
\numberwithin{equation}{section}

\newtheorem{theorem}{Theorem}
\numberwithin{theorem}{section}

\newtheorem{lemma}[theorem]{Lemma}

\def\eref#1{(\ref{#1})}

\newcommand{\R}{\mathbb{R}}
\newcommand{\C}{\mathbb{C}}

\newcommand{\N}{\mathbb{N}}
\newcommand{\HH}{\mathbb{H}}
\def\H{\mathbb{H}}
\def\I{I}
\def\diam{\mathrm{diam}}
\def\length{\mathrm{length}}

\def \A {{A}}
\def \downto {\searrow}
\def \eps {\epsilon}
\def \tx {{\tilde \xi}}
\def \P {{\bf P}}
\def\Bb#1#2{{\def\md{\bigm| }#1\bigl[\,#2\,\bigr]}}
\def\BB#1#2{{\def\md{\Bigm| }#1\Bigl[\,#2\,\Bigr]}}

\def\Pb{\Bb\P}
\def\Eb{\Bb\E}

\def\EB{\BB\E}

\def \p {{\partial}}
\def \E {{\bf E}}

\def\SG/{$SLE_\slepar$}
\def \SLE/{$SLE_6$}

\def \proof {{ \medbreak \noindent {\bf Proof.} }}

\def\llambda{\lambda_6}

\def \SLG/ {{\SG}}

\def\L{{\mathfrak{L}}}

\def\K{{\mathfrak K}}
\def \halpha {{\alpha}}
\def \hbeta {{\beta}}
\def\W{X}
\def \tx {{\tilde \xi}}

\long\def\hide#1{}

\def\ev#1{{\cal {#1}}} %events
\def\slepar{\kappa}
\def\st{\,:\,}
\def\closure{\overline}

\def \bar {\overline}

\def\E{{\bf E}}

\def\tilde{\widetilde}
\def\mobtr/{M\"obius transformation}
\def\hat{\widehat}
\def \l {{\ell}}

\begin{document}
\maketitle

\begin {abstract}
This paper determines values of intersection exponents between 
packs of planar Brownian motions in the half-plane and 
in the plane that were not derived in our first two papers.    
For instance, it is proven that the exponent 
$\xi (3,3)$ describing the asymptotic decay of
the probability of non-intersection between two packs of three independent 
planar Brownian motions each is 
$(73-2 \sqrt {73 }  ) / 12$.
More generally, the values of 
$\xi (w_1, \ldots, w_k)$ and 
$\tx (w_1', \ldots, w_k')$ are determined
for all $ k \ge 2$, $w_1 , w_2\ge 1$, 
$w_3 , \ldots,w_k\in[0,\infty)$ and
all $w_1',\dots,w_k'\in[0,\infty)$.
The proof  relies on the results  derived  in 
our first two papers and applies the same general methods.
We first find the 
two-sided exponents for the stochastic 
Loewner evolution processes in a half-plane, 
from which the Brownian
intersection exponents are determined
via a universality argument.
\end {abstract} 

%\newpage
%\tableofcontents
%\newpage

\section {Introduction}

This paper is a follow-up to the papers~\cite {LSW1,LSW2},
in which the exact values of many of
the intersection 
exponents between planar Brownian motions were determined.
It is assumed that the reader is familiar with the
terminology and the results of \cite{LSW1,LSW2}, to which we also 
refer for  
background (in particular, the link with 
critical exponents for other models,
such as critical percolation or self-avoiding 
walks in the plane) 
 and a more complete bibliography.

Let us first very  briefly recall  the definition of these intersection 
exponents. Suppose that $k \ge 2$, $n_1, \ldots, n_k \ge 1$
are integers, and that $(B^{l,j})_{1 \le l \le k,  1 \le j \le n_l}$ is a 
collection of independent planar Brownian motions started from 
distinct points in a half-plane ${H}$.
Define the $k$ packs of Brownian motions
$\mathfrak B^l(t):= \bigcup_{j=1}^{n_l} B^{l,j}[0,t]$,
$l=1,2,\dots,k$.
Consider the following events:
\begin{align*}
{\cal E}(t) 
=
{\cal E}_{(n_1,\dots,n_l)}(t)
&
:= 
\bigcap _{1\le l < l' \le k}
\bigl\{
\mathfrak B^{l} (t) \cap\mathfrak B^{l'} (t) = \emptyset \bigr\}
,
\\
\tilde {\cal E}(t) 
=
\tilde {\cal E}_{(n_1,\dots,n_l)}(t)
&:=
{\cal E}_{(n_1,\dots,n_l)}(t) \cap
\bigcap_{l=1}^k \bigl\{ 
\mathfrak B^{l} (t)  \subset {H} \bigr\}.
\end{align*}
It is easy to see, using a subadditivity 
argument, that when $t \to \infty$, 
$$\Pb{  {\cal E} (t)}
\approx t^{-\xi /2 },\qquad 
\Pb{  \tilde {\cal E}(t)} \approx t^{-\tilde \xi /2},
$$
for some  $\xi = \xi ( n_1, \ldots, n_k)$ 
and 
$ \tilde \xi = \xi (n_1, \ldots,n_k)$, which are
called
the intersection exponents between $k$ packs 
of $(n_1, \ldots , n_k)$ Brownian motions in the 
plane and in the half-plane, respectively. 
 Here, $f\approx g$ means $\lim_{t\to\infty}\log f/\log g=1$.

There exists natural extensions of $\xi$ and $\tilde \xi$ 
to non-integer values of $n_1, \ldots, n_k.$
For instance, one can define the exponents 
$\xi (1, w)$ and $\tx (1, w)$
for all $w>0$ by the relations
\begin{align*}
\EB{ \Pb { {\cal E}_{(1,1)} (t) \md \mathfrak B^1 (t) }^w }
&
\approx t^{-\xi (1, w) /2 },
\\
\EB{ \Pb { \tilde {\cal E}_{(1,1)} (t) \md \mathfrak B^1 (t) }^w }
&
\approx t^{-\tx (1, w) /2 }.
\end{align*}
It is easy to see that these exponents $\tx (1,w)$
and $\xi (1,w)$ coincide with the previously defined exponents
when $w$ is a positive integer.

A second generalization are the two-sided 
exponents $\tx (w, 1, w)$.
One way to define them is as follows:
Suppose that $k=3$, $n_1=n_2=n_3=1$, that 
$H$ is the upper half-plane $\H=\{ x+ i y \st y >0 \}$,
and that for all $l \in\{ 1,2,3\}$, 
$B^{l,1} (0) =  e^{i l \pi / 4}$ and define
$ \hat {\mathfrak B}^l (t) 
= (0, e^{il\pi/4}] \cup B^{l,1} [0,t]$
and the event $\hat {\cal E}(t)$ 
that $\hat {\mathfrak B}^l (t)$ for $l=1,2,3$ 
are disjoint subsets of the half-plane $H$.
Loosely speaking, adding the segments $(0, e^{il\pi/4}]$
ensures that the three Brownian motions maintain their cyclic order
around zero.
Then $\tx (w,1,w)$ is defined for all $w>0$ by 
$$
\EB{ \Pb { \hat {\cal E} (t) | \mathfrak B^2 (t) }^w }
\approx t^{-\tx (w,1, w) /2 }.
$$
These exponents 
$\tx (w,1,w)$ coincide with the above definition when $w$
is an integer \cite {LW1}.

It has been shown in \cite {LW1} that there
exists a unique extension 
of $\xi$
and $\tilde \xi$ to non-integer values of 
$n_1, \ldots, n_k$ that is symmetric in its arguments
and satisfies the ``cascade relations''
(for all $1<j<k-1$)
\begin{equation}\begin{aligned}\label{cascade} \xi(n_1,\ldots,n_k) &=
\xi\bigl(n_1,
   \ldots,n_j, \tx(n_{j+1},\ldots,n_k)\bigr) , \\
 \tx(n_1,\ldots,n_k) &= \tx\bigl(n_1,\ldots,n_j,
   \tx(n_{j+1},\ldots,n_k)\bigr) . 
\end{aligned}\end{equation}
The extension of $\tx$ is valid for all positive $n_1,
\ldots,n_k$ while the extension of $\xi$ also requires
that at least two of the arguments are greater
or equal to 1.

A first consequence \cite {LW1}
of these cascade relations is that the 
value of all extended exponents $\xi$ and $\tx$ 
(and in particular their 
values when $w_1, \ldots, w_k$
are positive integers) can be 
expressed in terms of the functions
$w \mapsto \xi (1, 1, w)$, $w \mapsto \tx (1, w)$ 
and $w \mapsto \tx ( w,1,w)$ defined for all $w >0$.

A second consequence \cite {LW1}
is that the (extended) full-plane exponents are expressible as
a function of the half-plane exponents
\begin {equation}
\label {xieta}
\xi(w_1,w_2,\dots,w_k) 
= \eta\bigl( \tilde\xi(w_1,w_2,\dots,w_k)\bigr),
\end {equation}
provided that $w_1,w_2\ge1$; however, the function $\eta$
was not  determined in~\cite{LW1}.

Set 
$$ 
U (x) = \sqrt {x + (1/24)} - \sqrt {1/24} .
$$
In \cite {LSW1}, we determined the function
$w \mapsto \tx(1/3,w)$, and using the cascade relations~\eref{cascade} 
concluded
that 
\begin{equation}\label{egen}
\tilde\xi(w_1,w_2,\dots,w_k)
= U^{-1}\Bigl( U(w_1)+U(w_2)+\cdots+ U(w_k)
\Bigr) \end{equation}
holds for all $k \ge 2$,
all $w_1,w_2,\dots w_{k-1}\in \bigl\{p(p+1)/6\st p\in\N\bigr\}$,
and all $w_k>0$. 
The equation~\eref{egen} is expands to
\begin {equation*}
\tx (w_1, \ldots , w_k )
=
\frac {
\left( 
\sqrt { 24 w_1 + 1} + \sqrt { 24 w_2 + 1 }
+ \cdots + \sqrt { 24 w_k +1 } 
- (k-1) \right)^2  - 1  }{ 24 } \,.
\end {equation*}
In \cite {LSW2}, we showed that $\xi (1,1) =5/4$,
determined the function $ w \mapsto \xi(1,1,1, w)$,
and concluded from the cascade relations and~\eref{egen}  that
\begin {equation}
\label {eta}
\forall x\ge7\qquad\eta(x)=\frac{\bigl(\sqrt{24x+1}-1\bigr)^2-4}{48}.
\end {equation} 
Combined with~\eref {egen} and~\eref{xieta}, this gives
the value of 
$\xi (w_1, \ldots, w_k)$ for a large collection of $w_1, \ldots, w_k$, but
not for all of them.

In the present paper, we will prove the following results.

\begin {theorem}
\label{tegen}
The identity \eref{egen} holds for all $k \ge 2$
and for all $w_1,\dots,w_k\ge 0$.
\end{theorem}

\begin {theorem} 
\label {tegen2}
The identity \eref {eta} holds for all $x \ge \tx (1,1)= 10/3$,
so that for all $k \ge 2$, $w_1, w_2 \ge 1$ and 
$w_3, \ldots,w_k \ge 0$,
\begin {eqnarray*}
\xi (w_1, \ldots, w_k ) 
&= &
\eta \circ U^{-1} \bigl( U(w_1) + \ldots + U (w_k)\bigr)
\\
&=&
\frac{ \left( \sqrt {24 w_1 + 1} + \cdots + \sqrt {24 w_k +1 }
- k \right)^2 - 4 
}{48
}\,. 
\end {eqnarray*}
\end {theorem}

These two theorems determine almost all the 
Brownian intersection exponents.
Those which they do not give are $\xi(n,w)$ where
$n\in\N_+$ and $ w\in(0,1)$. 
 It seems that the universality argument,
which is used to translate information about \SLE/ exponents to
Brownian exponents, cannot be extended to this range.
Therefore, the techniques of~\cite{LSW1,LSW2} and the present paper
do not suffice.

To complete the picture, in the forthcoming paper~\cite {LSWan}, we
determine these remaining
exponents by analytic continuation.  There, it will be shown
that for integers
$n \ge 1$, the mapping $w \mapsto \xi (n, w)$ 
is real-analytic 
in $(0, \infty)$. Combining this 
with the above theorems gives the 
value of $\xi (n,w)$ for all positive $w$
(i.e., removing the $w \ge 1$ condition), and gives the 
value of the disconnection exponents
$\xi(n,0):= \lim_{w \searrow 0} \xi (n, w)$ for all
$n\in\N_+$.  
That is,~\eref{xieta} also holds when $k=2$, $w_1\in\N_+$ and
$0\le w_2<1$ and~\eref{eta}
holds for all $x\ge 1$.  This will conclude the determination of
all the two-dimensional Brownian
 intersection exponents that have been
defined.

The proof of Theorem \ref {tegen} is very similar to the 
proofs we used to derive (\ref {egen}) and (\ref {eta})
in \cite {LSW1,LSW2}.
A crucial role is played by the stochastic Loewner evolution
process with parameter $6$ (\SLE/) introduced in \cite {S}. 
In the present paper, first the
two-sided exponents associated to \SG/ in a half-plane
are computed.
Then, via a universality argument, the values
of the Brownian half-plane exponents $\tx (1, w_1, 1, w_2)$
are deduced.
This leads directly to Theorem \ref {tegen} 
via the cascade relations satisfied by $\tx$.
Theorem \ref {tegen2} also immediately  follows by using  
 the results  
derived in \cite {LSW2}. 

Following is a rough and somewhat imprecise comparison 
of the approach used in the present paper in relation to
those of~\cite{LSW1} and~\cite{LSW2}.
In \cite{LSW1}, we have studied the expectation of the derivative
to any power $w\ge 0$ of a suitably normalized conformal map
onto the complement of a chordal \SG/ process crossing a rectangle
from left to right.
The expectation was determined precisely.  Its rate of decay as
a function of the width corresponds to the exponent
$\tx(1/3,w)$.  The reason that $1/3$ appears as the first argument,
rather than $1$, is that the \SG/ process was permitted to touch
one horizontal edge. 

In the present paper, we study the
expectation of an expression of the form
$f'(x_1)^{w_1}f'(x_2)^{w_2}$, where $f$ is a suitably
normalized conformal map.
The points $x_1$ and $x_2$ at which these derivatives
are computed are on the two sides of the \SG/ process,
hence the name of the paper.  The explicit formula for
the expectation is not calculated, however, the decay rate as
a function of the size of the \SG/
is determined, which suffices.  The decay rate
corresponds to the exponent $\tx(w_1,1,w_2)$.
The calculation of the decay rate is via an eigenvalue
computation, as in~\cite{LSW2}. 

In~\cite{LSW2}, the decay rate of expectation of a single derivative
raised to an arbitrary power $w_1>0$ was calculated for radial \SG/.

\section {Notations and terminology}

The present paper builds on the results of our previous papers
\cite{LSW1,LSW2},
and it will be
assumed that the reader is familiar with the terminology and
tools used in these papers. In particular, we refer to these two papers  
for definitions and properties of 
chordal and radial \SLE/, Brownian excursions in a domain, 
and their relation to Brownian intersection exponents. 

Let $f$ and $g$ be functions, and let $l\in\R$ or $l=\infty$.
Say that $ f(x) \sim g(x) $ when $x \to l$,
if $f(x)/ g(x)\to 1$.
Write $f(x) \approx g(x)$, if $\log f(x) / \log g(x)\to 1$,
and write $ f(x) \asymp g(x)$,
if $f(x) / g(x)$ is bounded above and below by 
positive finite constants when $x$ is sufficiently close to $l$.

For convenience, just as in \cite {LW1,LW2,LSW1,LSW2}, we will
use $\pi$-extremal distance, which is defined as $\pi$ times the 
usual extremal distance or extremal length in a domain. 
The $\pi$-extremal distance in a domain $D$ between two sets $A,A'$
will be denoted $\l(A,A';D)$.
For more information on extremal length, as well as other basic 
tools from complex analysis that we shall use (Koebe
$1/4$ Theorem, Schwarz Lemma), 
 see, for instance, \cite {A}.

\section{Derivative \SG/ exponents}
\label{sder}

Let $x\in(0,1)$, let $\slepar>0$,
let $K_t$ be the hulls of chordal \SG/ in $\HH$,
from $x$ to $\infty$, and let
$g_t:\HH\setminus K_t\to\HH$ be the
conformal maps normalized by the
hydrodynamic normalization
$\lim_{z\to\infty}g_t(z)-z =0$.
In other words, for all $z \in \HH$, $g_t (z)$ 
is the solution of the ordinary differential equation 
$$
\partial_t  g_t (z) = \frac { 2 } { g_t (z) -  W_t}
,\qquad  g_0 (z) = z, $$
where $t \mapsto W_{t / \slepar}$ is a standard
 real-valued Brownian motion started from 
$W_0 = x$.
The set $\HH \setminus K_{t_0}$ consists of all $z \in \HH$ such that 
$t \mapsto g_t (z)$ is well-defined at least up to time $t_0$.
Then, 
$(K_t, t \ge 0)$ is an increasing family of subsets of $\HH$:
For more details on the definition of $(K_t, t \ge 0)$, 
some of its properties such as scaling, conformal 
invariance,  see \cite {S,LSW1}.

Let 
$$T:=\inf\bigl\{t > 0 \st 
\{0,1\}\cap\closure{K_t}\neq\emptyset\bigr\}
$$
denote the first time at which the \SG/ swallows $0$ or $1$,
and $T:=\infty$ if no such time exists.
For all $t < T$,
let  
$$
f_t(z):= \frac{g_t(z)-g_t(0)}{g_t(1)-g_t(0)}
$$
be  the conformal map from 
$\HH \setminus K_t$ 
onto the upper half-plane  
such that $f_t (0) = 0, f_t (1) = 1 $ and $f_t (\infty) = \infty$.
Note that 
$ f_t' ( \infty ) = \bigl( g_t (1) - g_t (0)\bigr )^{-1}$
is decreasing and continuous 
in $t$ for $t < T$. Let $S:= - \lim_{t \nearrow T} 
\log f_t' (\infty)$.
Perform a time-change as follows: For all 
$s \in [0, S)$, define
$$
t(s) := \inf \bigl \{ t \in [0,t)
 \st f_t' ( \infty) \le   e^{-s}\bigr \}
$$
and the inverse map
$$ s(t) := -\log  f_t'  (\infty )$$
for all $t \in [0,T)$. 
For all $t<T$ define also 
$$
Y_{s(t)} =  Z_t := \frac { W_t - g_t ( 0)}{g_t (1) - g_t (0) }.
$$
(This $Z_t$ was already used in~\cite {LSW1}.)
Loosely speaking, $Y_s$ and $Z_t$
correspond to the image (under $f_t$) of the point where $K_t$
grows at time $t$. 

For all $s<S$, also set
$$ 
 \halpha (s )
:= -\log f_{t(s)}' (0),\qquad 
\hbeta (s) = -\log f_{t(s)}' (1).
$$
For every $w_1, w_2>0$ and 
 every smooth function $F : [0,1] \to [0,1]$, let 
$$
h_F (x,s)
= h_F (x, s, w_1, w_2 ) 
:= \E_x \Bigl[1_{\{s < S\}} 
 F( Y_s)
\exp\bigl( - w_1 \halpha(s) - w_2\hbeta(s)\bigr)\Bigr],
$$
where $\E_x$ refers to expectation with respect
to the \SG/ started at $x$; that is, $W_0=Z_0= Y_0 =x$.
In particular, write $h_1$ in case $F$ is the constant function $1$.
That is,
$$
h_1 (x,s) =
\E_x
\Bigl[
1_{ \{t(s) < T \}}\, f_{t(s)}'(0)^{w_1}\, f_{t(s)}'(1)^{w_2} 
\Bigr] .
$$

\begin{theorem}\label{deriv}
For all  $w_1,w_2 > 0$ and $\slepar>0$, there exists some
$c>0$ such that for all $x\in(0,1)$ and all $s\ge 1$
$$
G(x)\exp(-\lambda s)
\le h_1(x,s)
\le c\, G(x)\exp(-\lambda s),
$$
where
\begin{align*}
&
\lambda =\lambda_\slepar(w_1,w_2)
:=
\frac{\Bigl(
\sqrt{(\slepar-4)^2+16\slepar w_1}
+
\sqrt{(\slepar-4)^2+16\slepar w_2}
+\slepar\Bigr)^2-(8-\slepar)^2}
{16\slepar}\,,
\\
&
G(x)
:= x^{a_1}(1-x)^{a_2}
\,,
\qquad\qquad
a_j
:= \frac{\slepar -4 + \sqrt{(4-\slepar)^2+16 w_j\slepar}}{2\slepar},
\qquad j=1,2.
\end{align*}
\end{theorem}

\noindent
{\bf Remark.}
It can be shown that this
theorem also holds when $w_1 =0$ and/or $w_2 = 0$, but this
will not be done here.

\noindent
\proof
This is a first eigenvalue 
computation, and the proof will follow quite closely the 
proof of Lemma 3.2 in \cite {LSW2} (which is the corresponding result 
for radial \SG/, but with the derivative computed at only one point).
A simple computation (using the definitions of $\alpha$, $\beta$,
$Y$, $g_t$ and $s(t)$) shows that for all $s < S$,
\begin {equation}
\label {Y}
dY_s 
= \sqrt { \frac { \slepar Y_s (1-Y_s) }{2}}\: dB_s
+ (1 - 2 Y_s )\: ds ,
\end {equation}
where $B$ is a standard Brownian motion.
Note also that $S$ is the first time at which $Y$ hits
$\{0,1\}$, (unless $S=T=\infty$), and that 
\begin{equation}\label{ab}
\partial_s \halpha (s) = 1 / Y_s 
,\qquad
\partial_s \hbeta (s) = 1 / (1- Y_s) .
\end{equation}
We first use this to prove
that 
\begin{equation}\label{formu}
h_G(x,s) = \exp(-\lambda s) G(x)\,.
\end{equation}
Let $\W=[0,1]\times[0,\infty)$.
Set $h=h_G$ and let $\hat h(x,s) =\exp(-\lambda s) G(x)$.
Observe that 
$$
Q_s:=h\bigl(Y_s,s_0 - s \bigr)
\exp\bigl(-w_1\alpha(s)-w_2\beta(s)\bigr)
$$
is a local
martingale on $s \le  s_0$. 
(For this, the choice of $G$ is not important.)
Moreover, $h$ is smooth in $(0,1)\times(0,\infty)$.
Consequently, the $ds$ term in It\^o's formula 
for $dQ_s$ must vanish; that is,
\begin{equation}
\label{newimproved}
\p_s h =
(1-2x) \,\partial_x h  + (1-x)x\frac{\slepar}4 \,\partial_x^2 h
  -\Bigl( \frac{w_1}{x}  + \frac{w_2}{1-x}\Bigr) h\,.
\end{equation}
It is immediate to verify that $\hat h$
satisfies this differential equation in the interior of $\W$.
It is also clear that $\hat h= h$ on $\partial \W$.

In a moment, we shall see that $h$ is continuous in $\W$.
Assuming this for now, an easy application of
the maximum principle gives $h=\hat h$.
Indeed, let $\eps>0$, and suppose that there is some point $(x_0,s_0)$
with $h-\hat h\ge\eps$.
Among all such points, choose one with $s_0$ minimal.
Since $h=\hat h$ on $\partial \W$, $(x_0,s_0)$ must be in the interior.
By minimality of $s_0$, it follows that
$\p_s h(x_0,s_0)-\p_s\hat h(x_0,s_0)\ge 0$
and that $h(x,s_0)-\hat h(x,s_0)$
has a local maximum at $x_0$.  From the latter fact,
we may deduce that $\p_x (h-\hat h)=0$
and $\p_x^2 (h-\hat h)\le 0$ at $(x_0,s_0)$.
However, these facts put together contradict~\eref{newimproved},
and we may conclude that $h\le \hat h+\eps$.
The same argument shows that $h\ge\hat h-\eps$.
Since $\eps>0$ was arbitrary, it follows that $\hat h=h$.

To establish \eref{formu}, it therefore 
remains to prove the continuity of
$h$.
Suppose that $Y$ starts at $Y_0=x$ where $0<x < 2^{-n_0-2} \min \{s_0, 1\}$,
for some constants $s_0>0$ and $n_0\in\N_+$.
Define the stopping times $\nu_0 = 0$ and 
for all $n \ge 0$,
$$
\nu_{n+1} 
:= \inf \bigl\{ s > \nu_n \st s = \nu_n + Y_{\nu_n} 
\hbox { or } 
\left| Y_s  - Y_{\nu_n} \right|  \ge Y_{\nu_n} /2 \bigr\} .
$$
Note that for all $n \le  n_0-1$,
$ 0<  Y_{\nu_n} \le  2^{n} x $, 
$ \nu_{n+1} \le \sum_{j=0}^{n} Y_{\nu_j} \le 2^{n+1} x$,
so that $\nu_{n_0} \le s_0 $.
Let $\ev R_n$ denote the event
$$ 
\ev R_n 
:=
\{ \nu_{n} = \nu_{n-1} + Y_{\nu_n} \}.
$$
Let ${\cal F}_n$ denote
the $\sigma$-field generated by the events
$\ev R_1,\dots,\ev R_n$.
There is a $c>0$ such that for
$Y_s<1/2$, the diffusion term in~\eref{Y} is bounded below by
$c\sqrt {Y_s}$, 
and the drift term is bounded by $1$.
Hence,
it is not difficult
(for instance, using Girsanov's 
formula and Doob's inequality) to see that for all $n \le n_0$,
the conditional probability 
$\Pb{  \ev R_n \md {\cal F}_{n-1} }$
 is bounded below by a positive
constant, which does not depend on $x$ and $n$. On the event $\ev R_n$,  
we have $\halpha(\nu_n) - \halpha ( \nu_{n-1} ) \ge  2/3$,
by~\eref{ab}.
It follows easily that  
$\alpha (s_0)$ tends in probability to $\infty$ when $x \searrow 0$,
and therefore (since $w_1>0$)
 that $h$ tends to zero as $x\searrow0$ (uniformly 
for $s \ge s_0$).  A similar argument
shows that $h\to0$ as $x\nearrow 1$.  It is 
easy  to verify that for any $\eps>0$,
$h(x,s)\to G(x)$
as $s\to0$ uniformly with respect to $x \in (\eps, 1- \eps)$.
It is  
also easy 
to check that $h(x,s) \to 0$ when 
$(x,s) \to (0,0)$ or $(x, s) \to (1, 0)$
(note that $G(0)= G(1) =0$).
This shows that $h$ is continuous in $X$ 
and concludes the proof of \eref{formu}.

Since $G\le 1$, it is clear that for all $s>0$ and $x \in (0,1)$,
\begin {equation}
\label {easy}
h_1 (x,s) \ge h_G (x,s) = e^{-\lambda s} G(x) .
\end {equation}
It  remains to prove
that 
$
\inf_{x\in(0,1)}\inf _{s\ge 1} h_G(x,s)/h_1(x,s)>0\,.
$
The Markov property at time $s-1$ 
shows that it suffices to establish
this for $s=1$. Since both $h_1 ( \cdot, 1)$ 
and $h_G (\cdot, 1)$ are positive and continuous 
on $(0,1)$, it suffices to prove this when $x$ is close
to 0 and close to 1. By symmetry, it is enough to treat
the case where $x$ is close to $0$.
Now assume that
$0 < x \leq 1/2$.  For every positive integer
$n$, let 
$$
r_n:=\inf\Bigl\{
h_G(x,s)/h_1(x,s)\st
 x\in[4^{-n} , 1/2],\,s\in[1 - 2^{-n} ,1]\Bigr\}.
$$
 Assume $4^{-n} \leq x  <  4^{-n+1}$, $1 - 2^{-n}
\leq s \leq 1$, and let
\[  \tau := \inf\bigl\{ s\st Y_s \in \{0,4^{-n+1}\} \bigr\} . \]  
Note that~\eref{ab} gives
\[ \BB{\E_x}{e^{-w_1 \alpha(2^{-n})} 1_{\{ \tau > 2^{-n} \} }}
   \leq e^{-w_1 2^{n-2} } . \]
Since $h_1(x,s) \geq h_G(x,s) \geq c x^{a_1} $ for some
constant $c>0$ and all $(x,s)$ as chosen above, it follows that
\[ \Bb{\E_x}{e^{-w_1 \alpha(s) - w_2 \beta(s)} 1_{\{ s<S,\,\tau \leq 2^{-n}
    \}} }
    \geq (1 - \epsilon_n) \Bb{\E_x}{e^{-w_1 \alpha(s) - w_2 \beta(s)}
1_{\{ s<S \} }  } ,  \]
where $\epsilon_n:=c^{-1}\,4^{na_1}e^{-w_1 2^{n-2}}$.
However, since
$s - \tau \geq 1 - 2^{-n+1}$ on the event $\{\tau \leq 2^{-n}\}$,
and since $\{ s< S ,\, \tau \le 2^{-n} \}
\subset \{ Y_\tau  = 4^{-n+1} \}$, 
the strong Markov property gives
\begin {align*}
h_G(x,s)&\ge
{\BB {\E_x}{e^{-w_1 \alpha(s) - w_2 \beta(s)} \,G(Y_s) \,
1_{\{s<S,\, \tau \leq 2^{-n}\}} } }
\\&
=
\BB{\E_x}{e^{-w_1\alpha(\tau)-w_2\beta(\tau)} h_G(4^{-n+1},s-\tau)
1_{\{\tau\leq 2^{-n}, Y_\tau = 4^{-n+1}\}}}
\\&
\ge
r_{n-1}\,
\BB{\E_x}{e^{-w_1\alpha(\tau)-w_2\beta(\tau)} h_1(4^{-n+1},s-\tau)
1_{\{\tau\leq 2^{-n}, Y_\tau = 4^{-n+1}\}}}
\\&
=
     r_{n-1} \,
\BB{\E_x}{e^{-w_1 \alpha(s) - w_2 \beta(s)} 1_{\{  s < S ,\,\tau \leq 2^{-n}
    \}} } 
\\
& \geq 
 r_{n-1} \,(1- \epsilon_n)\,
   \BB{\E_x}{e^{-w_1 \alpha(s) - w_2 \beta(s)} 1_{ \{ s < S \} } } 
\\&
=
 r_{n-1} \,(1- \epsilon_n)\,
h_1(x,s)\,.
\end {align*}
That is, $r_n \geq (1 -\epsilon_n) r_{n-1}$.  Since $\sum_n\eps_n<\infty$,
this gives $\inf_n r_n > 0$, which completes the proof.
\qed

\section{Extremal distance exponents}

In the previous section, we derived estimates concerning 
the joint law
of $\log f'(0)$ and $\log f'(1)$ at the first time at which
$ f_t'(\infty) = e^{-s}$. We now use this result to 
obtain information concerning the law of the extremal distances 
at the first time at which \SG/ reaches distance $R$.
More precisely, let $R\ge 1$, and let
$V_R$ denote the half disk
$$
V_R:= \bigl\{z\in\H\st |z-1/2|<R\bigr\}\,.
$$
Let $A_R$ denote the semi-circle
$\H\cap\p V_R$.  Let $a\in(0,1)$, and
consider chordal \SG/ in $\H$ from $a$ to $\infty$.
Let 
$$
\tau=\tau_R:=
\inf\bigl\{t\st \closure K_t\cap A_R \not= \emptyset\bigr\}\,,
$$
and set 
$$
\K=\K_R:= \bigcup_{t<\tau} K_t\,.
$$
As before, let $T$ be the first time that the \SG/ swallows $0$
or $1$.
Let $\I_1(t):=[0,a]\setminus \closure K_t$ and
$\I_2(t):=[a,1]\setminus \closure K_t$.
On the event $\tau<T$,
let $\L_1(R):=\l(\I_1 (\tau) ,A_R;\H \setminus \K )$ denote the
$\pi$-extremal distance from 
$\I_1(\tau)$ to $A_R$ in $\H \setminus \K$
(or in $V_R \setminus \K$ since they are equal), and
let
$\L_2(R):=\l(\I_2 (\tau) ,A_R;\H \setminus \K)$.
 Let
\[  H(a,R) = \E_a\Bigl[1_{\{ \tau < T \} }
\exp\bigl(-w_1\L_1(R)-w_2\L_2(R)\bigr)\Bigr]
.\]

\begin{theorem}\label{el}
Let $\slepar>0$, 
$w_1,w_2 > 0$, and let $\lambda = \lambda_{\slepar}
( w_1, w_2)$ be as in Theorem~\ref {deriv}.
There is a constant $c= c ( \slepar,  w_1, w_2)$ such that for all 
$R>2$,
$$
\forall a\in(0,1)\quad
H(a,r)
\le c\, R^{-\lambda} . 
$$
On the other hand, for all 
  $a_0\in (0,1/2)$, there
is a $c'=c'(\slepar,w_1,w_2,a_0) >0$ such that for all $R>2$
$$
\forall a\in[a_0,1-a_0]\quad
H(a,r)
\ge c' \, R^{-\lambda}.
$$
\end{theorem}

\proof
We use the notation of Section~\ref{sder}.
Using scaling invariance and a monotonicity argument, it is easy to 
see that
for all $R > 2$ 
and $a \in [a_0, 1-a_0]$,
$$ 
H \bigl( 1/2 , (R+1) / a_0 \bigr) \le H (a, R)
\le H \bigl(1/2 , (R-1) / 2 \bigr) ,
$$
and hence it suffices to show that 
$ H (1/2 , R ) \asymp R^{- \lambda}$.

We now assume that $a = 1/2$. Let 
$$
\sigma_R:=t({\log R})=
\sup \bigl\{t < T \st  f_t'(\infty) > 1/R\bigr\}\,.
$$
For $t < T$, let $\tilde K_t$ be the union of $[0,1]$
with $K_t$ and with the reflection of $K_t$
about the real axis.
Observe that $f_t$ extends conformally
to a map $f_t:\C\setminus \tilde K_t\to\C\setminus[0,1]$.
Therefore, the Koebe 1/4 Theorem and the Schwarz Lemma
together imply that there is a constant $c_1>0$
such that for all $t < T$,
$$
c_1^{-1}\diam(\tilde K_t)\le f_t'(\infty)^{-1}  \le c_1 \,\diam(\tilde
K_t)\,.
$$
Since $R\le\diam(\tilde K_\tau)\le 2R$,
this gives
\begin{equation}\label{sigmatau}
\tau_{R/2c_1} \le \sigma_R \le \tau_{c_1 R} 
\end{equation}
for all $R \ge c_1$ such that $ \sigma_R < T $.
The Koebe 1/4 Theorem and the Schwarz Lemma then also 
show that if $\tau = \tau_R < T$,
\begin {equation}
\label {KS}
\diam\bigl(f_\tau(A_R)\bigr) \asymp 1 
,\qquad 
\diam\bigl(f_\tau(A_{2R})\bigr)\asymp 1 
.\end {equation}

For all $t < T$, let $L_1(t)$ be the length
 of the image of $I_1 (t)$ under
$f_t$,
and let $L_2(t)$ be the length of the image of
$I_2 (t)$ under $f_t$.
Recall that $Z_t=Y_{s(t)}$, and
note that $\p_t \log f'_t(x)$ is monotone decreasing in $x$ when
$f_t(x)\le Z_t$, and monotone increasing for $f_t(x)\ge Z_t$,
because
\begin{align*}
&\p_t  \log f_t' (x) 
+ \p_t  \log\bigl( g_t(1) - g_t(0)\bigr)
= \p_t \log g_t'(x) = \frac{\p_x\p_t g_t(x)}{g_t'(x)}
\\ &
{} \qquad
= \frac { -2 }{(g_t(x) - W_t)^2}  
= \frac { -2  }{(f_t (x) - Z_t)^2 (g_t(1) - g_t(0))^2} 
\,.
\end{align*}
Therefore, $f'_t(x)\le f'_t(0)$ for $x\in \I_1 (t)$
and $f'_t(x)\le f'_t(1)$ for $x\in \I_2 (t)$.
Consequently, $L_1(t)\le f_t'(0)/2$ and $L_2(t)\le f_t'(1)/2$.

Since the $\pi$-extremal distance between $A_R$ and $A_{2R}$ in $\H$
is $1$, it follows from (\ref {KS})
 that the Euclidean distance between
$f_\tau(A_R)$ and $f_\tau(A_{2R})$ is bounded away
from zero.  
Therefore
the Euclidean distance between $f_\tau(A_{2R})$ and
$f_\tau\bigl(\I_j (\tau )\bigr)$ is  bounded away from
 zero  for $j=1$ and $j=2$.
It is also bounded from infinity since the diameter 
of $f_\tau ( A_{2R})$ is bounded.
This shows that 
(provided $\tau < T$)
\begin{equation}\label{a2}
\exp\bigl(-\l(A_{2R},\I_1(\tau);\H\setminus \K_R)\bigr)\asymp 
 \length\bigl(f_\tau(\I_1 (\tau))\bigr)\le
 f'_\tau(0)/2,
\end{equation}
and similarly for $\I_2$.
Note that
$$
\L_1(R)
\le\l\bigl(A_{2R},\I_1(\tau_R);\H\setminus\K_R\bigr)
\le\l\bigl(A_{2R},\I_1(\tau_{2R});\H\setminus\K_{2R}\bigr)
= \L_1(2R).
$$
Hence, 
$$
\E\Bigl[
1_{\{ \tau_{2R} < T\} }
\exp\bigl(-w_1\L_1(2R)-w_2\L_2(2R)\bigr)\Bigr]\le c_2 R^{-\lambda}
$$
follows from Theorem~\ref{deriv}, \eref{sigmatau},
and \eref{a2}.

For the other direction, since $f'_t$ is monotone decreasing on
$(-\infty,1/2)\setminus K_t$, if $\tau < T $, then
$$
\length \bigl( f_\tau( [ -1, 1/2 ]\setminus\K  )\bigr ) 
\ge f_\tau' (0) 
$$
and
$$
\length \bigl( f_\tau ([1/2, 2]\setminus\K)\bigr)
\ge f_\tau' (1) 
.$$
Hence,
when $\tau <T$,
\begin {align*}
f_\tau' (0)
& \le
\length \bigl( f_\tau( [ -1, 1/2 ]\setminus\K  )\bigr ) 
\\
&
 \asymp 
 \exp\bigl(-\ell ( A_{2R} , [-1, 1/2] ; \H \setminus \K_R )\bigr)\\
& \le  
 \exp\bigl(-\ell ( A_{R} , [-1, 1/2] ; \H \setminus \K_R )\bigr),
\end {align*}
and similarly for $f_\tau' (1)$. 
Combining this with scaling invariance of \SG/ and 
Theorem \ref {deriv} then readily shows that 
$$
\E\Bigl[
1_{\{ \tau < T\} }
\exp\bigl(-w_1\L_1(R)-w_2\L_2(R)\bigr)\Bigr]\ge c_3 \,  R^{-\lambda},
$$
and completes the proof of Theorem \ref {el}.
\qed

\section{The universality argument}

Let $\mu_R$ denote the Brownian excursion measure
in the domain $V_R$, and let $B$ denote an excursion.
(See \cite {LW2, LSW1, LSW2} for the definition 
of the excursion measures on simply connected domain
 and the link with the Brownian
intersection exponents). 
Let $\ev Q_B$ denote the
event that the initial point $B(0)$ of
$B$ is in $(0,1/2)$, and the terminal point in
$A_R$.  On this event, let $\L$ be the
$\pi$-extremal distance from $[0,B(0)]$ to
$A_R$ in $V_R\setminus B$, and
let $\L_B$  be the $\pi$-extremal
distance from $[B(0),1]$ to $A_R$
in $V_R\setminus B$.
Then 
for all $w,w'>0$, when $R \to \infty$,
$$
\int_{\ev Q_B}
\exp ( - w \L - w' \L_B ) \: d\mu_R (B)
\approx R^{- \tx (w, 1, w')}.
$$

Let
$\phi=\phi_B$ be the conformal map
from the component $\W=\W_B$ of $V_R\setminus B$
whose boundary contains $[B(0),1]$
to a semi-disk $V_{\tilde R(B)}$ such that
$\phi$ takes $\p \W\cap [B(0),1]$ onto
$[0,1]$ and takes $\p \W\cap A_R$
onto $A_{\tilde R(B)}$.
Set $\tilde\L_B:=\log \tilde R(B)$.
Note that when $R \to \infty$,
$$ \tilde \L_B = \L_B +  O(1).
$$
Hence,
for all $w,w'>0$, when $R \to \infty$,
$$
\int_{\ev Q_B}
\exp ( - w \L - w' \tilde \L_B ) \: d\mu_R (B)
\approx R^{- \tx (w, 1, w')}.
$$

We will need a lemma saying that we can   
restrict ourselves to the case where $B(0) < 1/2$
and the  
conformal map $\phi$ does not 
push $1/2$ too close to $0$.
More precisely, let $\ev H$ denote the event that $\ev Q_B$ holds and
$\phi(1/2)\in[1/20,19/20]$.

\begin{lemma}\label{behaved}
For all $w,w'\ge 0$, as $R\to\infty$,
\begin{equation}\label{ll}
\int_{\ev H} \exp(-w \L-w'\L_B)\,d\mu_R(B)
\approx R^{- \tx (w,1,w')}\,.
\end{equation}
\end{lemma}

\proof
Let $M:=\{z\in\H\st |z-1|\le 9/10\}$,
and let $\ev M$ be the event $B\cap M=\emptyset$.
We first show that $\ev Q_B \cap\ev M\subset\ev H$.
Indeed, extend $\phi$ to $\{\bar z\st z\in \W_B\}$,
by Schwarz reflection.  Since $\phi(1)=1$ and
$\phi(\W_B)\supset \W_B$, it follows
from the Schwarz Lemma that $\phi(x)\le x$ for all
$x\in [B(0),1]$.  In particular
$\phi(1/2)\le 1/2$. 
Let $\psi$ be the conformal map
from the disk $\{z\st |z-1|< 9/10\}$ onto
$\C\setminus (-\infty,0]$ such that $\psi(1)=1$
and $\psi'(1)>0$.  The Schwarz Lemma also
shows that $\phi(x)\ge\psi(x)$ for all $x\in[1/10,1]$.
Since $\psi(z)=(10z-1)^2/(10z-19)^2$,
it follows that
$\phi(1/2)\ge\psi(1/2)>1/20$.  This proves
$\ev Q_B\cap\ev M\subset\ev H$.

On the event $\ev Q_B \cap \ev M$, let
$\L_B':=\l\bigl([B(0),1/10],A_R;V_R\setminus(B\cup M)\bigr)$
be the $\pi$-extremal distance
from $[B(0),1/10]$ to $A_R$ in $V_R\setminus(B\cup M)$.
Let $L'$ be the $\pi$-extremal distance from
$[0,1/10]$ to $A_R$ in $V_R\setminus M$.
It is clear that $\log R\le L'\le\log R+O(1)$.
Consequently, by the restriction property and
conformal invariance for
Brownian excursions, it follows that
\begin{equation}\label{upr}
\int_{\ev M\cap\ev Q_B} \exp(-w\L-w'\L_B')\,d\mu(B)
\approx R^{-\tx(w,1,w')}\,.
\end{equation}
Observe that $\L_B\le\L_B'$ and that
\begin{equation}\label{lowr}
\int_{\ev Q_B} \exp(-w\L-w'\L_B)\,d\mu(B)
\approx R^{-\tx(w,1,w')}\,.
\end{equation}
Since the left hand side of~\eref{ll}
is between the left hand sides of~\eref{upr} and~\eref{lowr},
the lemma follows.
\qed
\bigskip

\noindent
{\bf Proof of Theorem \ref {tegen}.}
Let $\mu_R$ denote the Brownian excursion measure in the
domain $V_R$, and let $B$ be an excursion.
Let $\P_R$ denote the law of \SLE/ in $V_R$ started from $1/2$,
and let $\K$ be as in the previous section.
Let $\ev Q_B$ be the event
that the initial point $B(0)$ of $B$ is in $[0,1/2]$
and the terminal point is in $A_R$,
let $\ev Q_\K$ be the event that $\K\subset\H\cup[0,1]$,
and let $\ev Q$ be the event $\ev Q_\K\cap\ev Q_B\cap \{\K\cap
B=\emptyset\}$.
On $\ev Q$
let
\begin{align*}
\L
&
:=\l\bigl([0,B(0)],A_R;V_R\setminus B\bigr),
\\
\L_B
&
:=\l\bigl([B(0),1],A_R;V_R\setminus B\bigr),
\\
\L_\K
&
:=\l\bigl([0,1/2],A_R;V_R\setminus\K\bigr),
\\
\L'
&
:=\l\bigl( [ B(0), 1/2 ] ,A_R; V_R \setminus (\K \cup B)\bigr),
\\
\L''
&
:=\l\bigl([1/2,1],A_R;V_R\setminus\K\bigr).
\end{align*}
We determine the
asymptotics as $R\to\infty$ of
$\Eb{1_{\ev Q}\exp(-w\L-w'\L'- w'' \L'')}$ in two different ways.

Given $B\in\ev Q_B$, we may map  $\W_B$ 
onto $V_{\exp(\tilde \L_B)}$
by $\phi$.
By conformal invariance of \SLE/, the restriction property 
of 
\SLE/ (see \cite {LSW1}) and Theorem~\ref{el}, it follows that 
$$
\Eb{\exp(-w'\L' - w'' \L'') \md B} 
\le c \exp\bigl(-\llambda(w', w'')\tilde \L_B\bigr)\,.
$$
Hence,
\begin{align*}
&{\int_{\ev Q_B} \int_{\ev Q_\K} 1_{\ev Q}
\exp(-w\L-w'\L' - w'' \L'')\,d\P_R(\K)\,d\mu_R(B)
}
\\ &
\quad \le  c \int_{\ev Q_B}
\exp\bigl(-w\L-\llambda(w', w'')\L_B\bigr)\,d\mu_R(B)
\\
&
\quad \approx   R^{-\tx(w,1,\llambda(w', w''))}\,.
\end{align*}
For the other direction, by Theorem \ref {el} 
and by Lemma~\ref{behaved}, we have
\begin{eqnarray*}
\lefteqn {
\int_{\ev Q_B} \int_{\ev Q_\K} 1_{\ev Q}
\exp(-w\L-w'\L'- w'' \L'')\,d\P_R(\K)\,d\mu_R(B)
}\\
& \ge&
\int_{\ev H} \int_{\ev Q_\K} 1_{\ev Q}
\exp(-w\L-w'\L'- w'' \L'')\,d\P_R(\K)\,d\mu_R(B)
\\
& \ge 
& c'\, \int_{\ev H}
\exp(-w\L)
\exp \bigl(-\llambda(w',w'')\tilde \L_B\bigr)\,d\mu_R(B)
\\
&\approx &
R^{-\tx(w,1,\llambda(w'w''))}\,.
\end{eqnarray*}
We may therefore conclude that
\begin{equation}\label{c1}
\int_{\ev Q_B} \int_{\ev Q_\K} 1_{\ev Q}
\exp(-w\L-w'\L'-w''\L'')\,d\P_R(\K)\,d\mu_R(B)
\approx R^{-\tx(w,1,\llambda(w', w''))}\,.
\end{equation}

On the other hand, by conformal invariance and the 
restriction property of the Brownian excursions,
given $\K\in\ev Q_\K$, we have
$$
\int 1_\ev Q \exp(-w\L - w' \L')\,d\mu_R(B)\approx
\exp\bigl(-\tx(w,1, w')\L_\K\bigr).
$$
Consequently, Theorem~\ref{el} gives
\begin{eqnarray*}
\lefteqn {\int_{\ev Q_\K}
\int_{\ev Q_B}
1_{\ev Q}
\exp(-w\L-w'\L'- w'' \L'')
\,d\mu_R(B)
\,d\P_R(\K)
}
\\
&
\approx
&\int_{\ev Q_\K}
\exp\bigl(-\tx(w,1, w')\L_\K-w''\L''\bigr)
\,d\P_R(\K)
\\&
\approx
&R^{-\llambda(\tx(w,1, w'),w'')}.
\end{eqnarray*}
Comparing with~\eref{c1} gives
\begin{equation}\label{cook}
\llambda\bigl(\tx(w,1, w'),w''\bigr) =
\tx\bigl(w,1,\llambda(w',w'')\bigr)\,.
\end{equation}

Define $y(w') := \lim_{w\searrow0}\llambda (w, w')$.
First let $w \downto 0$ and $w' \downto 0$ in (\ref {cook}). Recall that 
$(w, w') \mapsto \tx (w,1,w')$ is continuous at $(0,0)$
(see e.g., \cite {LW1}), 
so that for all $w'' \ge 0$ 
\begin {equation}
\label {der1}
\llambda ( 1, w'')=
\tx \bigl(1, y(w'')\bigr) ,
\end {equation}
which shows that $\tx (1, v) 
= \llambda \bigl(1, y^{-1} (v)\bigr)$ in the case where 
 $v \ge 1$ (we also derived this result in \cite {LSW1}).

Now $w \downto 0$ and $w'' \downto 0$ in~\eref{cook} gives 
$$
 y\bigl( \tx (1, w')\bigr)
=
\tx \bigl(1, y (w')\bigr) .
$$
Combining this with (\ref {der1})
and the explicit expression for $\llambda$ shows that 
for all $v>0$
\begin {equation}
\label {der2}
\tx (1, v) =
y^{-1} \bigl(\llambda (1,v)\bigr)
=
y(v).\end {equation}

Finally, letting $w' \downto 0$ in~(\ref {cook}) shows that 
$
 \tx \bigl(w, 1, y( w'')\bigr)
= \llambda \bigl( \tx (w, 1) , w''\bigr)$,
which gives
$$
\tx \bigl(w,1,\tx(1,w'')\bigr)
 = \llambda \bigl( y(w), w''\bigr).$$
The cascade relations~\eref{cascade} and (\ref {der2})
applied to the left hand side imply
$$
\tx (w,w'') = y^{-1}\circ y^{-1} \circ 
  \llambda \bigl( y(w), w''\bigr).$$
Via further applications of the cascade relations,
this
leads to the explicit expression for all
$\tx (w_1, \ldots, w_k)$.
\qed

\medbreak

\noindent 
{\bf Proof of Theorem \ref {tegen2}.}
By (\ref {xieta}) and Theorem \ref {tegen}, it suffices to 
derive the value of $\xi (1, w, 1, w)$ for all $w >0$.
This is a simple combination of Theorem \ref {el}, the relation between 
radial and chordal \SLE/ (see \cite {LSW2}), and the computation of 
exponents for radial \SLE/ (see \cite {LSW2}). The proof is 
essentially the same as  in the final section of \cite {LSW2}.
One has to consider (for small $r>0$) a Brownian excursion 
$B$ in the annulus $\A_r = \{ z \st r < |z| < 1 \}$ and an independent
radial \SLE/ started at $1$ (growing towards 0) stopped
when it hits the circle of radius $r$, and the event $\ev C$ 
that they both cross the annulus without intersecting each other. Define
$\L$ and $\L'$ to be the two $\pi$-extremal
distances between the two circles in each of 
the two connected components of $ \A_r \setminus ( B \cup \K)$ that 
cross the annulus. The result is derived by estimating the integral of 
$\exp ( - w \L - w \L')$ in two ways. First,  fixing $B$ and   
applying Theorem \ref {el}  and  \cite[Lemma 5.5]{LSW2} 
gives the 
exponent $\xi ( 1, \llambda (w,w))$;  this is equal to $\xi (1,
\tx (w,1,w) ) = \xi (1,1,w,w)$ by 
Theorem \ref {tegen}.
On the other hand, if we first fix $\K$ and use the radial \SLE/ exponents 
derived in \cite {LSW2} and the value of $\tx (w,1,w)$ (from Theorem 
\ref {tegen}), we can compute explicitly the exponent.
Since this is almost word for word the same argument as in 
the final section of \cite {LSW2}, we safely leave the details 
to the reader.
\qed
\medskip

The fact that \eref {eta} is 
valid for all $x \ge \tx (1,1)$
is also an immediate corollary of~\eref{eta}
and the analyticity result from~\cite{LSWan}.
Consequently, Theorem~\ref{tegen2} also follows from~\cite{LSWan}
and Theorem~\ref{tegen}.

\bigskip

\filbreak
\begingroup
\small 
\parindent=0pt 

\ifhyper\def\email#1{\href{mailto:#1}{\texttt{#1}}}\else
\def\email#1{\texttt{#1}}\fi
%\noindent\hbox{
%\noindent
\vtop{
\hsize=2.3in
Greg Lawler\\
Department of Mathematics\\
Box 90320\\
Duke University\\
Durham NC 27708-0320, USA\\
\email{jose@math.duke.edu}
}
\bigskip
\vtop{
\hsize=2.3in
Oded Schramm\\
Microsoft Corporation,\\
One Microsoft Way,\\
Redmond, WA 98052; USA\\
\email{schramm@Microsoft.com}
}
\bigskip
\vtop{
\hsize=2.3in
Wendelin Werner\\
D\'epartement de Math\'ematiques\\
B\^at. 425\\
Universit\'e Paris-Sud\\
91405 ORSAY cedex, France\\
\email{wendelin.werner@math.u-psud.fr}
}
%}
\endgroup

\filbreak

\end {document}